\documentclass[12pt,a4paper]{article}
\usepackage{amssymb}

\usepackage{graphicx,amssymb,amsfonts,epsfig,amsthm,a4,amsmath,url}
\usepackage[latin1]{inputenc}

\newtheorem{Thm}{Theorem}[section]
\newtheorem{Prop}[Thm]{Proposition}
\newtheorem{Lem}[Thm]{Lemma}
\newtheorem{Cor}[Thm]{Corollary}
\newtheorem{Ex}[Thm]{Example}
\theoremstyle{definition}
\newtheorem{Def}[Thm]{Definition}

\newtheorem{Rem}[Thm]{Remark}

\newcommand{\Z}{\mathbf{Z}}

\newcommand{\CC}{\mathcal{C}}
\newcommand{\N}{\mathbf{N}}
\newcommand{\R}{\mathbf{R}}
\newcommand{\C}{\mathbf{C}}
\newcommand{\Q}{\mathbf{Q}}

\newcommand{\T}{\mathbf{T}}

\title{The Chabauty space of $\Q_p^\times$}
\author{Antoine Bourquin and Alain Valette}
\date{\today}

\begin{document}

\baselineskip=16pt

\maketitle

\begin{abstract} 
\textit{Let $\CC(G)$ denote the Chabauty space of closed subgroups of the locally compact group $G$. In this paper, we first prove that $\CC (\Q_p^\times)$ is a proper compactification of $\N$, identified with the set $N$ of open subgroups with finite index. Then we  identify the space $\CC(\Q_p^\times) \smallsetminus N$ up to homeomorphism: e.g. for $p=2$, it is the Cantor space on which 2 copies of $\overline{\N}$ (the 1-point compactification of $\N$) are glued.}
\end{abstract}

\section{Introduction}

In 1950, Chabauty \cite{Cha} introduced a topology on the set $\mathcal{F}(X)$ of closed subsets of a locally compact space $X$, turning $\mathcal{F}(X)$ into a compact space, see Definition 1 below; for $X$ discrete, this is nothing but the product topology on $2^X$. When $G$ is a locally compact group, the set $\CC(G)$ of closed subgroups of $G$ is a closed subset of $\mathcal{F}(G)$, so $\CC(G)$ is a compact set canonically associated with $G$: we call it the {\it Chabauty space} of $G$.

\begin{Def} For a locally compact space $X$, the Chabauty topology on $\mathcal{F}(X)$ has as open sets finite intersections of subsets of the form
$$\mathcal{O}_K=\{F\in\mathcal{F}(X): F\cap K=\emptyset\}$$
with $K$ compact in $X$, and
$$\mathcal{O}'_U=\{F\in\mathcal{F}(X): F\cap U\neq\emptyset\}$$
with $U$ open in $X$.
\end{Def}

Let us give some examples of Chabauty spaces, first for additive groups of some locally compact fields:
\begin{Ex}\begin{enumerate}
\item For $G=\R$, the Chabauty space $\CC(\R)$ is homeomorphic to a closed interval, say $[0,+\infty]$, with the subgroup $\lambda\Z$ (with $\lambda >0$) being mapped to $\lambda$, the subgroup $\{0\}$ being mapped to $+\infty$, and the subgroup $\R$ being mapped to 0.
\item For $G=\C$, the situation is already much more subtle, and it was proved by Hubbard and Pourezza \cite{HP} that $\CC(\C)$ is homeomorphic to the 4-sphere.
\item For $G=\Q_p$, the field of $p$-adic numbers, every non-trivial closed subgroup is of the form $p^k\Z_p$ for some $k\in\Z$, so $\CC(\Q_p)$ is homeomorphic to the 2-point compactification of $\Z$, namely $\Z\cup\{\pm\infty\}$, with the subgroup $p^k\Z_p$ being mapped to $k\in\Z$, the subgroup $\{0\}$ being mapped to $+\infty$ and the subgroup $\Q_p$ being mapped to $-\infty$.
\end{enumerate}
\end{Ex}

Let us turn to multiplicative groups of some locally compact fields:
\begin{Ex}\begin{enumerate}
\item For $G=\R^\times$, since $G\simeq\R\times\Z/2\Z$, we get 3 types of closed subgroups:
\begin{itemize}
\item Closed subgroups of $\R$, contributing a copy of $[0,+\infty]$.
\item Subgroups which are products $H\times \Z/2\Z$, with $H$ a closed subgroup of $\R$. They contribute another copy of $[0,+\infty]$, with origin the subgroup $\R\times\Z/2\Z$ and extremity the subgroup $\{0\} \times\Z/2\Z$.
\item Infinite cyclic subgroups which are not contained in $\R$; those are of the form $\langle(\lambda,1)\rangle$, for some $\lambda>0$. So they contribute a copy of $]0,+\infty[$.
\end{itemize}
As $\langle(\lambda,1)\rangle$ converges to $\R\times\Z/2\Z$ for $\lambda\rightarrow 0$, and to $\{(0,0)\}$ for $\lambda\rightarrow +\infty$, we see that subgroups of the 3rd type ``connect'' subgroups of the first and the second type, so that $\CC(\R^\times)$ is homeomorphic to a closed interval.
\item For $G=\C^\times\simeq \R\times\T$, the structure of $\CC(G)$ was elucidated by Haettel \cite{Hae}: it is path connected but not locally connected, and with uncountable fundamental group.
\end{enumerate}
\end{Ex}

In this paper, we deal with the multiplicative group $G=\Q_p^\times$ of the field $\Q_p$. We enjoy the general results by Y. Cornulier \cite{Cor} on $\CC(H)$ for $H$ a locally compact abelian group. Applied to $G=\Q_p^\times$, they yield that $\CC(G)$ is totally disconnected and uncountable (Theorems 1.5 and 1.6 in \cite{Cor}), that $\{1\}$ is not isolated in $\CC(G)$ while $\{G\}$ is isolated in $\CC(G)$ (Lemmas 4.1 and 4.2 in \cite{Cor}). More generally, a closed subgroup $H$ defines an isolated point in $\CC(G)$ if and only if $H$ is open with finite index in $G$ (Theorem 1.7 in \cite{Cor}).

To describe our main result, recall that a compact metric space $Y$ is a {\it proper compactification of} $\N$ if $Y$ contains an open, countable, dense, discrete subset $N$. It is known (see e.g. Propositions 2.1 and 2.3 in \cite{Tsan}) that every non-empty compact metric space can be written as  $Y\smallsetminus N$, with $Y$ and $N$ as above, in a unique way up to homeomorphism. 

We will use the following notations: $\overline{\N}=\N\cup\{\infty\}$ is the one-point compactification of $\N$, $C$ is the Cantor space, $[k]$ is the set $\{1,2,...,k\}$, and $d(n)$ is the number of divisors of $n$. 

\begin{Thm}\label{main} Let $p$ be a prime. For $G=\Q_p^\times$:
\begin{enumerate} 
\item $\CC(G)$ is a proper compactification of $\N$, viewed as the set $N$ of open subgroups of finite index.
\item For $p$ odd, the space $\CC(G)\smallsetminus N$ is homeomorphic to the space obtained by glueing $[d(p-1)]\times \overline{\N}$ on $C$, with the $d(p-1)$ accumulation points of $[d(p-1)]\times \overline{\N}$ being identified to $d(p-1)$ pairwise distinct points of $C$.
\item For $p=2$, the space $\CC(G)\smallsetminus N$ is homeomorphic to the space obtained by glueing $[2]\times \overline{\N}$ on $C$, with the 2 accumulation points of $[2]\times \overline{\N}$ being identified to 2 pairwise distinct points of $C$.
\end{enumerate}
\end{Thm}

In the above picture, the isolated points of $\CC(G)\smallsetminus N$ correspond to the closed infinite subgroups of $\Z_p^\times$, the invertible group of the ring $\Z_p$ of $p$-adic integers; and (for $p$ odd) the glueing points are the finite subgroups of $\Q_p^\times$, which are exactly the cyclic groups $C_d$ of $d$-roots of unity, for $d$ dividing $p-1$. It follows from the result that every non-isolated point of $\CC(G)\smallsetminus N$ is a condensation point, and the Cantor-Bendixson rank of  $\CC(G)\smallsetminus N$ is 1.

%\begin{Rem}
%For the case $p=2$, $\CC(G)$ is also a compactification of $\N$ and $\CC (G) \smallsetminus N$ is homeomorphic to space obtained by glueing $\{1,2\} \times \overline{\N}$ on $C$ the same way as for $p \neq 2$.
%\end{Rem}

To prove Theorem \ref{main}, we use the decomposition as topological groups $\Q_p^\times\simeq \Z_p\times C_{p-1}\times\Z$ (for $p$ odd), and $\Q_2^\times\simeq\Z_2\times C_2\times\Z$ (for all this, see section 3 in Chapter II of \cite{Ser}). For a locally compact abelian group $H$ with Pontryagin dual $\hat{H}=\hom(H,\T)$, the space $\CC(H)$ identifies with $\CC(\hat{H})$ for $H$: see Proposition \ref{orthogonal} below for the precise statement of this result of Cornulier \cite{Cor}. So, actually we work with the Pontryagin dual of $\Q_p^\times$, which for $p$ odd identifies with $C_{p^\infty}\times C_{p-1}\times\T$, where $C_{p^\infty}$ denotes the $p$-Pr\"ufer group.

At the end of the Introduction in \cite{Cor}, Cornulier mentions as a non-trivial question to determine the homeomorphism type of $\CC(G)$, when $G$ is a totally disconnected locally compact abelian group that fits into a short exact sequence 
$$0\rightarrow A\times P\rightarrow G\rightarrow D\rightarrow 0,$$
where $A$ is finitely generated abelian, $D$ is Artinian (i.e. it has no infinite decreasing chain of subgroups), and $P$ is a finite direct product of finite groups and groups isomorphic to $\Z_p$ (for various $p$'s, multiplicities allowed). Since $\Q_p^\times$ is of that form, Theorem \ref{main} and its extension Theorem \ref{technicalmain} can be seen as a contribution to Cornulier's question.

The structure of the paper is as follows: in section 2, we gather some preliminary material about the Chabauty space $\CC(G)$. In section 3, we determine the Chabauty space of $C(C_n\times\Z)$: the Chabauty space of any finitely generated abelian group has been determined in Theorem C of \cite{CGP2}, but we give a direct proof for completeness in the case of $C_n\times\Z$; this will be needed in the proof of our main result. Although we do not need it, we see that the determination of the Chabauty space of $\Z^2$ follows easily from the one of $C_n\times\Z$. Finally, in section 4 we prove Theorem \ref{main} along the lines sketched above.

\medskip
\noindent
{\it Acknowledgements:} We thank Y. Cornulier and P. de la Harpe for useful comments on a first draft of this paper.

\section{Preliminaries}
\subsection{General locally compact groups}

In this subsection $G$ will denote a locally compact second countable group. First, we recall a result on convergence in $\CC (G)$ that we will use throughout the paper without further reference (see \cite{Pau} for a proof) :

\begin{Prop}
A sequence $(H_n)_{n>0}$ in $ \CC (G)$ converges to $H$ if and only if 
\begin{enumerate}
\item $\forall h \in H$, $\exists h_n \in H_n$ such that $h = \lim h_n$.
\item For every sequence $(h_{n_k})_{k>0}$ in $G$ with $h_{n_k} \in H_{n_k}$, $n_{k+1} >n_k$, converging to some $h \in G$, we have $h \in H$.
\hfill$\square$
\end{enumerate}
\end{Prop}

The following is lemma 1.3(1) in \cite{CGP1}.

\begin{Lem}\label{quotient} Let $N\triangleleft G$ be a closed normal subgroup, let $p:G\rightarrow G/N$ denote the quotient map. The map
$$p^*:\CC(G/N)\rightarrow \CC(G):H\mapsto p^{-1}(H)$$
is a homeomorphism onto its image.
\hfill$\square$
\end{Lem}

The proof of the next lemma is an easy exercise.

\begin{Lem}\label{compactnormal} Let $N\triangleleft G$ be a compact normal subgroup. If $(H_k)_{k>0}$ converges to $H$ in $\CC(G/N)$, then $(H_kN)_{k>0}$ converges to $HN$ in $\CC(G)$.
\hfill$\square$
\end{Lem}

If $L$ is a closed subgroup of $G$, the map $\CC(G)\rightarrow\CC(L):H\mapsto H\cap L$ is in general {\it not} continuous (take e.g. $G=\R^2$, and $L$ a 1-dimensional subspace). We single out a case where it is continuous.

\begin{Lem}\label{continuous} Let $G=D\times L$ be a locally compact group, with $D$ discrete. Identify $L$ with $\{1\}\times L$. The map $\CC(G)\rightarrow\CC(L):H\mapsto H\cap L$ is continuous.
\end{Lem}

{\bf Proof:} Assume that the sequence $(H_n)_{n>0}$ converges to $H$ in $\CC(G)$, let us check that $(H_n\cap L)_{n>0}$ converges to $H\cap L$ in $\CC(L)$.
\begin{itemize}
\item If $(1,l)$ is in $H\cap L$, we can write $(1,l)=\lim_{n\rightarrow\infty}(d_n,l_n)$, with $(d_n,l_n)\in H_n$ for every $n>0$. Since $D$ is discrete, we have $d_n=1$ for $n$ large enough, so that we have also $(1,l)=\lim_{n\rightarrow\infty}(1,l_n)$, with $(1,l_n)\in H_n\cap L$ for every $n$.
\item If a sequence $((1,l_{n_k}))_{k>0}$, with $(1,l_{n_k})\in H_{n_k}\cap L$, converges to $(1,l)\in L$, then by assumption we have $(1,l)\in H$, i.e. $(1,l)\in H\cap L$.\hfill$\square$
\end{itemize}

If $G$ is a locally compact {\it abelian} group, recall that $\hat{G}=\hom(G,\T)$ denote its Pontryagin dual. For $H$ a closed subgroup of $G$, denote by $H^\perp$ the {\it orthogonal} of $H$:
$$H^\perp=\{\chi\in\hat{G}:\chi|_H=1\}.$$
The following result is due to Cornulier (Theorem 1.1 in \cite{Cor}):

\begin{Prop}\label{orthogonal} The orthogonal map $\CC(G)\rightarrow \CC(\hat{G}):H\mapsto H^\perp$ is an inclusion-reversing homeomorphism.
\hfill$\square$
\end{Prop}

\subsection{Discrete groups}

The symbol $k\gg 0$ means ``{\it for $k$ large enough}''. In this section, $G$ denotes a discrete group. Let $(P_k)_{k>0}, P$ denote subsets of $G$. We say that $P$ is {\it locally contained} in $(P_k)_{k>0}$ if, for every finite subset $S\subset P$, we have $S \subset P_k$ for $k\gg 0$. In other words, $P$ is locally contained in $(P_k)_{k>0}$ if every finite subset of $P$ is contained in all but finitely many $P_k$'s.

\begin{Lem}\label{fglimit} Let $(H_k)_{k>0}, H$ be subgroups of $G$.
\begin{enumerate}
\item The sequence $(H_k)_{k>0}$ converges to $H$ in $\CC(G)$ if and only if $H$ is locally contained in $(H_k)_{k>0}$ and $G\smallsetminus H$ is locally contained in $(G\smallsetminus H_k)_{k>0}$.
\item If $(H_k)_{k>0}$ converges in $\CC(G)$ to $H$ and $H$ is finitely generated, then $H\subset H_k$ for $k\gg 0$.
\item The sequence $(H_k)_{k>0}$ converges to the trivial subgroup $\{e\}$ of $G$, if and only if for every finite subset $F\subset G\smallsetminus\{e\}$, we have $F\cap H_k=\emptyset$ for $k\gg 0$.
\end{enumerate}
\end{Lem}

{\bf Proof:} \begin{enumerate}
\item Assume that $(H_k)_{k>0}$ converges to $H$. If $S$ is a finite subset of $H$, then for any $g\in S$ we find a sequence $(h_k)_{k>0}$ in $G$, with $h_k\in H_k$ for every $k$, such that $g=\lim_{k\rightarrow\infty}h_k$. As $G$ is discrete, this means $g=h_k$ for $k\gg0$, i.e. $g\in H_k$ for $k\gg0$. As $S$ is finite, we have $S\subset H_k$ for $k\gg 0$. Similarly, if $T$ is finite and disjoint from $H$, we must show that $T$ is disjoint from all but finitely many $H_k$'s. Suppose not, i.e. $T$ intersects infinitely many $H_k$'s. As $T$ is finite, we find $g\in T$ and a subsequence $(k_i)_{i>0}$ such that $g\in H_{k_i}$ for every $i>0$. Setting $g_{k_i}:=g$, and writing $g=\lim_{i\rightarrow \infty} g_{k_i}$, because $(H_k)_{k>0}$ converges to $H$ this implies $g\in H$, which contradicts $T\cap H=\emptyset$.

Conversely, if $H$ is locally contained in $(H_k)_{k>0}$ and $G\smallsetminus H$ is locally contained in $(G\smallsetminus H_k)_{k>0}$, let us check that $(H_k)_{k>0}$ converges to $H$. For $g\in H$, set $h_k=g$ so that $h_k\in H_k$ for $k\gg 0$: so we write $g$ as a limit of elements in the $(H_k)'s$. Now, let $({k_i})_{i>0}$ be a subsequence and $h_{k_i}\in H_{k_i}$ such that $(h_{k_i})_{i>0}$ converges to $g\in G$. This means $g=h_{k_i}$ for $i\gg 0$, so that $g$ belongs to infinitely many $H_k$'s. As $G\smallsetminus H$ is locally contained in $(G\smallsetminus H_k)_{k>0}$, this forces $g\in H$, hence $\lim_{k\rightarrow\infty} H_k=H$.

%$\{g_1,...,g_m\}$ be a finite generating set of $H$. So for every $i=1,...,m$ there exists a sequence $(h_{i,k})_{k>0}$, with $h_{i,k}\in H_k$, converging to $g_i$. As $G$ is discrete this means $g_i=h_{i,k}$ for $k\gg 0$. So $H\subset H_k$ for $k\gg 0$. 
\item Apply the previous point to a finite generating $S$ of $H$: it is contained in $H_k$ for $k\gg 0$, so the same holds for $H$.
\item The condition is equivalent to $G\smallsetminus\{e\}$ being locally contained in $(G\smallsetminus H_k)_{k>0}$.

\hfill$\square$
\end{enumerate}

\begin{Lem}\label{findex} Let $G$ be a finitely generated group.
\begin{enumerate}
\item Any finite index subgroup defines an isolated point in $\CC(G)$.
\item The index of a subgroup is continuous on $\CC(G)$. More precisely, the map $\CC(G)\rightarrow\overline{\N}:H\mapsto [G:H]$ is continuous.
\end{enumerate}
\end{Lem}

{\bf Proof:} 
\begin{enumerate}
\item Let $H$ be a finite index subgroup of $G$ and let $(H_k)_{k>0}$ be a sequence in $\CC(G)$ converging to $H$. Let us show that eventually $H_k=H$. As $H$ is finitely generated, by Lemma \ref{fglimit} we already know that $H\subset H_k$ for $k\gg 0$. But there are finitely many distinct subgroups containing $H$, as $H$ has finite index. Passing to a subsequence we may assume that $H_k=K$ for $k\gg 0$. So the sequence $(H_k)_{k>0}$ converges both to $H$ and $K$, so $H=K$.
\item Let $(H_k)_{k>0} \subset \CC(G)$ be a sequence converging to $ H \in \CC(G)$; let us prove that $[G:H_k] \to [G:H]$. We have to consider two cases :
\begin{itemize}
\item Let assume that $[G:H] = n < \infty$. By the first part, $H$ is isolated in $\CC(G)$, so $H_k=H$ for $k\gg 0$, and the result follows.
%As a consequence, $H$ is finitely generated. By Lemma \ref{fglimit}, we know that $H \subset H_k$ for $k \gg 0$ and we can assume that $H \subset H_k$ for all $k>0$. Because of the hypothesis $[G:H]< \infty$, $H_k/H$ is finite. So we can assume that $H_k = L$ for all $k$. So the sequence $(H_k)_{k>0}$ converges both to $H$ and $L$ so $H=L$. We have that $[G:H_k] = [G:H]$.
\item Now assume that $[G:H] = \infty$.
If $[G:H_k] = \infty$ for $k \gg 0$, the result follows. Else, we may assume that $[G:H_k] = n_k < \infty $. If $n_k \not  \to  \infty$, passing to a subsequence, we may assume that $n_k =n$ for all $k$. As $G$ is finitely generated, the number of subgroups of index $n$ is finite. By the pigeonhole principle, we may assume that the sequence $(H_k)_{k>0}$ is constant. Hence, $[G:H]<\infty$ and this is a contradiction.
\hfill$\square$
\end{itemize}

\end{enumerate}

Let $IC(G)$ denote the set of infinite cyclic subgroups of a group $G$.

\begin{Lem}\label{infcyc} Let $G=\Z^d\oplus F$ be a finitely generated abelian group, with $F$ a finite abelian group, written additively. 
\begin{enumerate}
\item The closure of $IC(G)$ in $\CC(G)$ is $IC(G)\cup\{(0,0)\}$.
\item Let $(g_k)_{k>0}$ be a sequence in $G$. The sequence $\langle g_k\rangle_{k>0}$ converges to $\{(0,0)\}$ if and only if $\lim_{k\rightarrow\infty}g_k=\infty$ in $G$.
\item The topology induced by $\CC(G)$ on $IC(G)$ is discrete.
\end{enumerate}
\end{Lem}

{\bf Proof:} \begin{enumerate}
\item Assume that the sequence $\langle g_k\rangle_{k>0}$ in $\CC(G)$ converges to a subgroup $H$. We know that every subgroup of $G$ is finitely generated and by Lemma \ref{fglimit}, $H \subset \langle g_k \rangle$ for $k \gg 0$. So $H$ is a subgroup of the cyclic infinite group $\langle g_k\rangle $. So $H$ is either trivial or infinite cyclic. 
\item This follows immediately from the last part of lemma \ref{fglimit}.
%If $(g_k)_{k>0}$ converges to infinity, then clearly $(<g_k>)_{k>0}$ converges to $\{(0,0)\}$. For the converse, assume that $(g_k)_{k>0}$ does not converge to infinity; passing to subsequences, we may first assume that the $g_k$'s stay in some bounded set, second (as bounded sets in $G$ are finite) that the $g_k$'s are constant, say equal to some $g\in G$ of infinite order. Then the sequence $(<g_k>)_{k>0}$ does not converge to $\{(0,0)\}$.
\item Let $(\langle g_k\rangle )_{k>0} \subset IC(G)$ be a sequence converging to $\langle g\rangle \in IC(G)$. Let us show that $\langle g_k\rangle =\langle g\rangle$ for $k \gg 0$. By Lemma \ref{fglimit}, we may assume that $\langle g\rangle \subset \langle g_k\rangle$ for $k>0$. So, there exists non-zero integers $m_k \in \Z$ such that $g = g_k m_k$. 
We may assume that $m_k >0$ up to replacing $g_k$ by $-g_k$. If $m_k \to \infty$, then $g$ is divisible by infinitely many integers, a contradiction. So the $m_k$'s are bounded and, passing to a subsequence, we may assume that $m_k  = m$ for all $k>0$, i.e. $g=mg_k$. Since the equation $g=mx$ has finitely many solutions in $G$, by the pigeonhole principle we may assume that the sequence $(g_k)_{k>0}$ is constant. Because of the convergence of $(\langle g_k\rangle)_{k>0}$ to $\langle g\rangle$, this implies $g_k\in \langle g\rangle$, i.e. $m=1$, and $g_k=g$.\hfill$\square$

%So, the element $h =g_k$ is an element of $<g>$ and that's a contradiction.
\end{enumerate}

\section{The case $\Z\times\Z/n\Z$ and $\Z^2$}

\subsection{The group $\Z\times\Z/n\Z$}

In this subsection, we set $G=\Z\times\Z/n\Z$. Recalling that the rank of a finitely generated group is the minimal number of generators, by viewing $G$ as a quotient of $\Z^2$ we see that every subgroup of $G$ has rank at most 2. Now $G$ has three types of subgroups:
\begin{itemize}
\item Finite subgroups, i.e. subgroups of $\{0\}\times \Z/n\Z$: there are exactly $d(n)$ of them, where $d(n)$ is the number of divisors of $n$;
\item Infinite cyclic subgroups;
\item Infinite subgroups of rank 2.
\end{itemize}

We describe the structure of infinite subgroups more precisely.

\begin{Lem}\label{structure} Let $H$ be an infinite subgroup of $G$. Set $F:=H\cap(\{0\}\times \Z/n\Z)$. There exists $(a,b)\in H$, with $a>0$, such that $H=F\oplus \langle(a,b)\rangle$. In particular, $H$ is cyclic if $F$ is trivial and $H$ has rank 2 if $F$ is non-trivial. 
\end{Lem}

{\bf Proof:} Let $p_1:G\rightarrow\Z$ be the projection onto the first factor. As $p_1(H)$ is infinite, it is of the form $p_1(H)=a\Z$ for some $a>0$. Let $(a,b)\in H$ be any element such that $p_1(a,b)=a$. Then clearly $F \cap\langle (a,b)\rangle =\{(0,0)\}$, and for $(x,y)\in H$, writing $x=ma$ for some integer $m$, we get $(x,y)=m(a,b)+(0,y-mb)$ so that $H=F \oplus \langle (a,b)\rangle$.
\hfill$\square$

\medskip
The group $G$ has the feature that every infinite subgroup has finite index, so defines an isolated point in $\CC(G)$, by Lemma \ref{findex}. By Lemma \ref{infcyc}, the only accumulation point of infinite cyclic groups is the trivial subgroup. So it remains to study accumulation points of rank 2 subgroups, which are necessarily finite subgroups. Every non-trivial finite subgroup $H$ of $\Z/n\Z$ is the limit of the sequence $(k\Z\times H)_{k>0}$ of rank 2 subgroups. The converse is provided by:

\begin{Prop}\label{Z/nZ} Let $m$ be a divisor of $n$. Consider a sequence $(H_k)_{k>0}$ of infinite subgroups of rank 2 of $G$. It converges in $\CC(G)$ to $\{0\}\times \langle m\rangle $ if and only if, for $k\gg 0$, there exists $g_k\in G$ with $\lim_{k\rightarrow\infty}g_k=\infty$, such that $H_k$ is generated by $(0,m)$ and $g_k$. 
\end{Prop}

{\bf Proof:} For the sufficient condition: if $H_k=\langle (0,m),g_k\rangle=\langle (0,m)\rangle \oplus\langle g_k\rangle $, as $\lim_{k\rightarrow\infty}\langle g_k\rangle =\{(0,0)\}$ by Lemma \ref{infcyc}, we have $\lim_{k\rightarrow\infty} H_k=\langle(0,m)\rangle$ by Lemma \ref{compactnormal}. 
\\For the necessary condition: assume $(H_k)_{k>0}$ converges to $\{0\}\times \langle m\rangle$. By Lemma \ref{structure}, we have $H_k=F_k\oplus \langle g_k\rangle$ with $F_k=H_k\cap(\{0\}\times \Z/n\Z)$ and $g_k=(a_k,b_k)$ with $a_k>0$. Because of the assumed convergence, we have $F_k=\langle(0,m)\rangle $ for $k\gg 0$, and $\lim_{k\rightarrow\infty} g_k=\infty$.
\hfill$\square$

\medskip
Recall that we denote by $[k]$ the set $\{1,2,...,k\}$, and by $d(n)$ the number of divisors of $n$. From Proposition \ref{Z/nZ}, we get immediately the following special case of Theorem C in \cite{CGP1}:

\begin{Cor}\label{C(Z/nZ)} The Chabauty space $\CC(\Z\times\Z/n\Z)$ is homeomorphic to $\overline{\N}\times [d(n)]$, the accumulation points corresponding to the subgroups $\{0\}\times\langle m\rangle$, with $m$ a divisor of $n$.
\hfill$\square$
\end{Cor}

\subsection{The group $\Z^2$}

The list of subgroups of $G=\Z^2$ is as follows:
\begin{itemize}
\item The trivial subgroup $\{(0,0)\}$;
\item Subgroups of rank 1, i.e. infinite cyclic subgroups;
\item Subgroups of rank 2, i.e. finite index subgroups in $G$ (which define isolated points in $\CC(G)$, by Lemma \ref{findex}).
\end{itemize}

In each infinite subgroup $H$ of $G$, pick a minimal vector $m_H$ (so $m_H$ has minimal norm among all non-zero vectors in $H$). From the 3rd part of Lemma \ref{fglimit}, we immediately get:

\begin{Prop} Let $(H_k)_{k>0}$ be a sequence of infinite subgroups of $G$. This sequence converges to $\{(0,0)\}$ in $\CC(G)$ if and only if $\lim_{k\rightarrow\infty}\|m_{H_k}\|=+\infty$.
\hfill$\square$
\end{Prop}

It remains to see how a rank 1 subgroup $\langle h\rangle$ can be a limit in $\CC(G)$.

\begin{Prop} A sequence $(H_k)_{k>0}$ in $\CC(G)$ converges to the rank 1 subgroup $H=\langle h\rangle$ if and only if, for $k\gg 0$, there exists $g_k\in G$ such that $H_k=\langle h,g_k\rangle $ and the sequence $(g_k+H)_{k>0}$ goes to infinity in $G/H$.
\end{Prop}

{\bf Proof:} Write $h=(a,b)$, let $n>0$ be the GCD of $a$ and $b$, set $p=(\frac{a}{n},\frac{b}{n})$, so that $p$ is a primitive vector proportional to $h$. Let $c,d\in\Z$ be integers such that $ad-bc=n$, so that, with $q=(c,d)$, the set $\{p,q\}$ is a basis of $G$, and every vector in $G$ may be written uniquely $\alpha p+\beta q$, for $\alpha,\beta\in\Z$. The map
$$\pi:\Z^2\rightarrow (\Z/n\Z)\times\Z: \alpha p+\beta q\mapsto (\alpha (\mbox{mod}\,n), \beta)$$
is then a surjective homomorphism with kernel $H$. 

If $(H_k)_{k>0}$ is a sequence of subgroups converging to $H$, we have $H\subset H_k$ for $k\gg 0$, by Lemma \ref{fglimit}. So to study convergence to $H$, we may as well assume that $H\subset H_k$ for every $k>0$. By Lemma \ref{quotient}, such a sequence $(H_k)_{k>0}$ converges to $H$ if and only if the sequence $(\pi(H_k))_{k>0}$ converges to the trivial subgroup in $\CC((\Z/n\Z)\times\Z)$. By Proposition \ref{Z/nZ}, this happens if and only if, for $k\gg 0$, the subgroup $\pi(H_k)$ is infinite cyclic, say $\pi(H_k)=\pi(\langle g_k\rangle)$ for some $g_k\in \pi^{-1}\pi(H_k)=H_k$, with the property that $\lim_{k\rightarrow \infty}\pi(g_k)=\infty$. This concludes the proof.
\hfill$\square$

\medskip
As a consequence, we get another special case of Theorem C in \cite{CGP1}:

\begin{Cor} The Chabauty space $\CC(\Z^2)$ is homeomorphic to $\overline{\N}^2$.
\hfill$\square$
\end{Cor}

\section{The proof of Theorem \ref{main}}

Let $p$ be a prime. Recall that $C_k$ denotes the cyclic group of order $k$ (viewed as the group of $k$-th roots of 1 in $\T$). It is classical (see \cite{Ser}) that 
$$\Q_p^\times \simeq \Z_p\times C_{p-1}\times\Z\;\;\mbox{($p$ odd)};$$
$$\Q_2^\times\simeq\Z_2\times C_2\times\Z.$$
Hence by Pontryagin duality
$$\widehat{\Q_p^\times}\simeq C_{p^\infty}\times C_{p-1}\times \T\;\;\mbox{($p$ odd)};$$
$$\widehat{\Q_2^\times}\simeq C_{2^\infty}\times C_2\times \T$$
where $C_{p^\infty}= \cup_{\ell=1}^\infty C_{p^\ell}$ denotes the Pr\" ufer $p$-group. By Proposition \ref{orthogonal}, $\Q_p^\times$ and $\widehat{\Q_p^\times}$ have canonically isomorphic Chabauty spaces. We choose to work with $\widehat{\Q_p^\times}$ as it is a 1-dimensional Lie group. More generally, for $k>0$ an integer, we set $G_{p,k}:=C_{p^\infty}\times C_k\times\T$ and we aim to determine $\CC(G_{p,k})$.

We will need some notation. We will denote by $\pi_1,\pi_2,\pi_3$ the projections of $G_{p,k}$ onto $C_{p^\infty}$ (resp. $C_k$, resp. $\T$). Set also $\pi=(\pi_1,\pi_2):G_{p,k}\rightarrow C_{p^\infty}\times C_k$.

We first give a list of closed subgroups of $G_{p,k}$: as it is a 1-dimensional Lie group, closed subgroups are either discrete, or 1-dimensional.

\begin{Lem}\label{1dim} Every closed 1-dimensional subgroup of $G_{p,k}$ is of the form $H_D^{(1)}:=D\times\T$, where $D$ is any subgroup of $C_{p^\infty}\times C_k$. The set of 1-dimensional subgroups is closed in $\CC(G_{p,k})$ and identifies with $\CC(C_{p^\infty}\times C_k)$ via $\pi^*$ (defined as in Lemma \ref{quotient}).
\end{Lem}

{\bf Proof:} A 1-dimensional subgroup of $G_{p,k}$ has the same connected component of identity as $G_{p,k}$, namely $\{1\}\times\{1\}\times\T$. The second statement follows immediately from Lemma \ref{quotient}.
\hfill$\square$

\medskip
\begin{Rem} Note that $\CC(C_{p^\infty}\times C_k)$ is homeomorphic to $\overline{\N}\times [d(k)]$, with accumulation points corresponding to the subgroups $C_{p^\infty}\times C_d$, for $d$ a divisor of $k$: this follows from Proposition \ref{Z/nZ} by dualizing.
\end{Rem}

We now turn to infinite discrete subgroups of $G_{p,k}$. For $F$ a finite subgroup of $C_k\times\T$, we denote by $q_F:G_{p,k}\rightarrow G_{p,k}/(\{1\}\times F)$ the quotient map. As we will need homomorphisms  $C_{p^\infty}\rightarrow (C_k\times\T)/F$ in Proposition \ref{graph} below, we start by describing those homomorphisms. We denote by $T^0$ the connected component of identity of the 1-dimensional compact abelian Lie group $(C_k\times\T)/F$. We choose some identification $T^0\simeq\T$, and we denote by $\iota: T^0\rightarrow (C_k\times\T)/F$ the inclusion. 

\begin{Lem}\label{homom} The map
$$\widehat{C_p^{\infty}}\rightarrow \hom(C_{p^\infty},(C_k\times\T)/F):f\mapsto \iota\circ f$$
is an isomorphism of compact groups.
\end{Lem}

{\bf Proof:} It is enough to see that any homomorphism $f:C_{p^\infty}\rightarrow (C_k\times\T)/F$ takes values in $T^0$. First observe that, being a 1-dimensional compact abelian group, $(C_k\times\T)/F$ is isomorphic to $A\times T^0$, for some finite abelian group $A$. Second, $C_{p^\infty}$ is a divisible group, and so is every homomorphic image of $C_{p^\infty}$. This applies in particular to the projection of $f(C_{p^\infty})$ to the first factor $A$. But a finite divisible group must be trivial, hence $f(C_{p^\infty})\subset T^0$.
\hfill$\square$

%\begin{Rem} We observe that $(C_k\times \T)/F$ is isomorphic to $C_d\times\T$, for some divisor $d$ of $k$: indeed $(C_k\times \T)/F$ is a 1-dimensional compact abelian Lie group, hence of the form $A\times\T$ for some finite abelian group $A$; looking at the groups of connected components we see that $A$ must be a quotient of $C_k$, i.e. of the form $C_d$ for some divisor $d$ of $k$. 
%\end{Rem}

\begin{Prop}\label{graph} For a finite subgroup $F$ of $C_k\times \T$ and a homomorphism $f:C_{p^\infty}\rightarrow (C_k\times\T)/F$, set  $H_{F,f}:=q_F^{-1}(Graph(f))$. Then $H_{F,f}$ is an infinite discrete subgroup of $G_{p,k}$ and every infinite discrete subgroup is of that form.
\end{Prop}

{\bf Proof:} It is clear that $H_{F,f}$ is an infinite discrete subgroup of $G_{p,k}$. Conversely, let $H$ be an infinite discrete subgroup. We first claim that $\pi_1(H)=C_{p^\infty}$: otherwise we would have that $\pi_1(H)=C_{p^\ell}$ for some $\ell>0$, so that $H$ appears as a discrete subgroup of the closed subgroup $\pi_1^{-1}(C_{p^\ell})$ of $G_{p,k}$: as this is a compact subgroup, this forces $H$ to be finite, a contradiction.

Set then $F:=\ker(\pi_1|_H)=H\cap(\{1\}\times C_k\times\T)$: as a discrete subgroup in a compact group, $F$ is finite. We consider two cases.

{\it Special case}: $F$ is trivial. Then $\pi_1|_H$ is injective, hence $f:=(\pi_2,\pi_3)\circ (\pi_1|_H)^{-1}$ is a homomorphism $C_{p^\infty}\rightarrow C_k\times\T$, and $H$ is exactly the graph of $f$.

{\it General case}: Let $F$ be arbitrary. We will reduce to our special case. We consider the infinite discrete subgroup $q_F(H)$ in $G_{p,k}/(\{1\}\times F)$. By the previous remark, for some divisor $d$ of $k$, we may identify $(C_k\times \T)/F$ with $C_d\times\T$, hence also $G_{p,k}/(\{1\}\times F)$ with $G_{p,d}$. By construction, the kernel $\ker(\pi_1|_{q_F(H)})$ is trivial, so we are back to the special case: therefore there exists a homomorphism $f:C_{p^\infty}\rightarrow (C_k\times \T)/F$ such that $q_F(H)=Graph(f)$. Then $H=q_F^{-1}(q_F(H))=H_{F,f}$ as claimed.
\hfill$\square$

\begin{Ex}\label{nonlift} (with $k=1$) Let $F=C_p$ viewed as a subgroup of $\T$. The quotient map $q_F$ identifies with the map $\T\rightarrow\T:z\mapsto z^p$. Let $\iota$ denote the inclusion of $C_{p^\infty}$ into $\T$. Then
$$H_{F,\iota}=\{(w,z)\in C_{p^\infty}\times\T: w=z^p\}.$$
Observe that $\iota$ does not lift, i.e. there is no $\tilde{\iota}:C_{p^\infty}\rightarrow\T$ such that $\iota=q_F\circ\tilde{\iota}$.
\end{Ex}

%\medskip
To summarize, from Proposition \ref{graph} and Lemma \ref{1dim}, we have now the complete list of closed subgroups of $G_{p,k}$:
\begin{itemize}
\item Finite groups;
\item Discrete infinite subgroups $H_{F,f}$, as described by Proposition \ref{graph};
\item 1-dimensional closed subgroups $H^{(1)}_D$, as described by Lemma \ref{1dim}.
\end{itemize}

The first part of Theorem \ref{main} follows from:

\begin{Prop}\label{compactification} $\CC(G_{p,k})$ is a proper compactification of $\N$, identified with the set $N$ of finite subgroups of $G_{p,k}$.
\end{Prop}

{\bf Proof:} By Theorem 1.7 in \cite{Cor}, every finite subgroup defines an isolated point in $\CC(G_{p,k})$. So $N$ is open and discrete in $\CC(G_{p,k})$, it remains to show that it is dense. This follows from:
\begin{itemize}
\item $H^{(1)}_D$ is the limit of the sequence $(D\times C_n)_{n>0}$ of finite subgroups.
\item $H_{F,f}$ is the limit of the sequence $(q_F^{-1}(Graph(f|_{C_{p^\ell}})))_{\ell>0}$, also consisting in finite subgroups.
\hfill$\square$
\end{itemize}

\medskip
Our aim is now to determine the homeomorphism type of $\CC(G_{p,k})\smallsetminus N$. Let $\CC^{(1)}$ be the set of one-dimensional subgroup, and let $\CC^{dis}$ be the set of infinite discrete subgroups, so that $\CC(G_{p,k})\smallsetminus N$ is the disjoint union of $\CC^{(1)}$ and $\CC^{dis}$. By lemma \ref{1dim} and the remark following it, $\CC^{(1)}$ is homeomorphic to $\overline{\N}\times[d(k)]$, the accumulation points corresponding to the subgroups $C_{p^\infty}\times C_d\times \T$, with $d$ a divisor of $k$. Note that the latter subgroup is also the limit of the sequence $(H_{C_d\times C_n,1_{d,n}})_{n>0}$ where $1_{d,n}$ denotes the trivial homomorphism $C_{p^\infty}\rightarrow (C_k\times\T)/(C_d\times C_n)$.

%To proceed, we observe that any homomorphism $f:C_{p^\infty}\rightarrow (C_k\times \T)/F$ as in Proposition \ref{graph}, may be viewed as a character of $C_{p^\infty}$. Indeed, identifying as above $(C_k\times \T)/F$ with $C_d\times\T$, for some divisor $d$ of $k$, we see that the projection of $f(C_{p^\infty})$ to $C_d$ is a $p$-group in $C_d$, hence it must be trivial because $p$ is prime to $k$. For each such finite subgroup $F$, we identify $Hom(C_{p^\infty},(C_k\times \T)/F)$ with $\widehat{C_{p^\infty}}$.

We denote by $\CC^{fin}(C_k\times\T)$ the set of finite subgroups of $C_k\times\T$. By the dual version of Corollary \ref{C(Z/nZ)}, this is a discrete set in $\CC(C_k\times\T)$.

\begin{Thm}\label{technicalmain}
\begin{enumerate}
\item The map 
$$\alpha: \CC^{fin}(C_k\times\T)\times\widehat{C_{p^\infty}}\rightarrow \CC^{dis}:(F,f)\mapsto H_{F,f}$$
is a homeomorphism.
\item The closure of $\CC^{dis}$ in $\CC(G_{p,k})$ is the union of $\CC^{dis}$ with the $d(k)$ subgroups $C_{p^\infty}\times C_d\times \T$, with $d$ a divisor of $k$.
\end{enumerate}
\end{Thm}

{\bf Proof:}\begin{enumerate}
\item \begin{itemize}
\item $\alpha$ is onto, by combining lemma \ref{homom} with Proposition \ref{graph}. 
\item To show that $\alpha$ is injective, we observe that $H_{F,f}$ determines both $F$ and $f$: first, $F$ is obtained as the intersection of $H_{F,f}$ with $\{1\}\times C_k\times \T$; second, $q_F(H_{F,f})$ is the graph of $f$, which of course determines $f$.
\item $\alpha$ is continuous. As $\CC^{fin}(C_k\times\T)$ is discrete, it is enough to show that, for each $F\in \CC^{fin}(C_k\times\T)$, the map $\widehat{C_{p^\infty}}\rightarrow \CC^{dis}:f\mapsto H_{F,f}$ is continuous.

Let fix $F \in \CC^{fin}(C_k\times\T)$ and let $(f_n)_{n>0} \subset \widehat{C_{p^\infty}}$ be a sequence converging to $f \in \widehat{C_{p^\infty}}$. We want to show that $H_{F,f_n} \to H_{F,f}$.
\begin{itemize}
\item Let $z \in H_{F,f}$. We have $q_F(z) = (a,f(a))$ for a certain $a \in C_{p^\infty}$. As $f_n \to f$, we have that $(a,f_n(a)) \to (a,f(a))$. So $q_F^{-1}(a,f_n(a)) \to q_F^{-1}(a,f(a))$ in the Chabauty topology of closed subsets. Moreover, $|F| = | q_F^{-1}(a,f_n(a))| = |q_F^{-1}(a,f(a))| < \infty$, so we can find $z_n \in H_{F,f_n}$ such that $z_n \to z$.

\item Let $z_{n_i} \in H_{F,f_{n_i}}$ such that $z_{n_i} \to z \in  G_{p,k}$ and let us show that $z \in H_{F_f}$. We have $q_F(z_{n_i}) = (a_{n_i},f(a_{n_i}))$ for a certain $a_{n_i} \in C_{p^\infty}$. By discreteness of $C_{p^\infty}$, convergence of $(z_{n_i})$ and continuity of $q_F$, we may assume that $a_{n_i} = a$ for all $i>0$. As $f_n \to f$, $q_F(z_{n_i}) \to (a,f(a))$. This implies that $z \in q_F^{-1} (a,f(a))$ and $z \in H_{F,f}$.
\end{itemize}

\item $\alpha$ is open. To show this, it is enough to show that, for every $F\in \CC^{fin}(C_k\times\T)$, the set $\{H_{F,f}:f\in\widehat{C_{p^\infty}}\}$ is open. This is in turn implied by the continuity of the map $\CC^{dis}\rightarrow\CC^{fin}(C_k\times\T):H_{F,f}\mapsto F$. Since $F$ can be expressed as $H_{F,f}\cap(\{1\}\times C_k\times \T)$, the result follows from lemma \ref{continuous}.
\end{itemize}

\item Let $(H_{F_n,f_n})_{n>0}$ be a sequence in $\CC^{dis}$ converging to $H\in\CC(G_{p,k})$. If the orders of the $F_n$'s remain bounded, passing to a sub-sequence we may assume that $F_n=F$ for $n\gg 0$. By compactness of $\widehat{C_{p^\infty}}$, again passing to a subsequence we may assume that $(f_n)_{n>0}$ converges to $f\in \widehat{C_{p^\infty}}$. Then $(H_{F_n,f_n})_{n>0}$ converges to $H_{F,f}$. Assume now that the orders of the $F_n$'s are unbounded. Passing to a subsequence, we may assume that $(F_n)_{n>0}$ converges in $\CC(C_k\times\T)$ to $C_d\times\T$, for some divisor $d$ of $k$. By the dual version of Proposition \ref{Z/nZ}, $F_n$ contains $C_d\times\{1\}$ for $n\gg 0$. For every $(\lambda_n,z_n)\in C_d\times \T$, the sequence of cosets $((\lambda_n,z_n)F_n)_{n>0}$ converges to $C_d\times\T$ in the Chabauty topology of closed subsets of $C_d\times\T$. So for $(w,\lambda,z)\in C_{p^\infty}\times C_d\times\T$, choose $(\lambda_n,z_n)$ in $q_{F_n}^{-1}(f_n(w))$, then choose $(\lambda'_n,z'_n)\in F_n$ such that $(\lambda_n\lambda'_n,z_nz'_n)$ converges to $(\lambda,z)$. So we have expressed $(w,\lambda,z)$ as the limit of the $(w,\lambda_n\lambda'_n,z_nz'_n)$'s in $H_{F_n,f_n}$. Conversely, if $(w_i,\lambda_i,z_i)\in H_{F_{n_i},f_{n_i}}$ converges to $(w,\lambda,z)\in G_{p,k}$, then $\lambda_i=\lambda$ for $i\gg 0$, hence $\lambda\in C_d$, and $(w,\lambda,z)\in C_{p^\infty}\times C_d\times\T$.
\hfill$\square$
\end{enumerate}

\medskip

The second part of Theorem \ref{main} is then a special case of :

\begin{Cor}\label{Coro_Q_p} $\CC(G_{p,k})\smallsetminus N$ is homeomorphic to the space obtained by glueing $[d(k)]\times \overline{\N}$ on the Cantor space $C$, with the $d(k)$ accumulation points of $[d(k)]\times \overline{\N}$ being identified to $d(k)$ pairwise distinct points of $C$.
\end{Cor}

{\bf Proof:}  By Theorem \ref{technicalmain}, the closure of $\CC^{dis}$ is a metrizable compact space which is totally disconnected and perfect (no isolated point), so it is homeomorphic to the Cantor space $C$. On the other hand we already observed after Proposition \ref{compactification} that $\CC^{(1)}$ is homeomorphic to $\overline{\N}\times[d(k)]$,  and the $d(k)$ accumulation points are identified with $d(k)$ pairwise distinct points of $C$.
\hfill$\square$

\bigskip
Authors'addresses:

\medskip
\noindent
Institut de Math\'ematiques\\
Universit\'e de Neuch\^atel\\
11 Rue Emile Argand\\
CH-2000 Neuch\^atel - Switzerland\\
antoine.bourquin@unine.ch\\
alain.valette@unine.ch

\end{document}